\documentclass[12pt]{article}
\usepackage{amsmath}
\usepackage{amsfonts,euscript,amssymb,amsmath,amsthm}
\usepackage{amstext}
\usepackage[english]{babel}

\newcommand{\bbC}{\mathbb{C}}

\newcommand{\Real}{\mathop{\rm Re}}

\newenvironment{pfl1}{\par\noindent\textbf{Proof of lemma 1.}}{\hfill$\square$}

\newtheorem{theorem}{Theorem}
\newtheorem{lemma}{Lemma}

\begin{document}
\large

\begin{center}
\section*{One metric result about analytic continuation of some
Dirichlet series.}

I.S. Rezvyakova
\footnote{This work was supported by grant RFBR 06-01-00518a and grant of the President of Russian Federation
MK-3370.2007.1}
\end{center}

We consider
$$
f(s,\theta) = \prod_{p} \left( 1-\frac{e^{2\pi i\theta p}}{p^{s}}\right)^{-1}, \quad \Real s >1, \quad
0 \le \theta \le 1,
$$
where the product is taken over all prime numbers, i.e.
$$
f(s, \theta) = \sum\limits_{n=1}^{+\infty} a_n(\theta) n^{-s},
$$
and for $n=p_1^{\alpha_1}\ldots p_k^{\alpha_k}$
$$
a_n(\theta) = e^{2\pi i(\alpha_1 p_1 + \ldots + \alpha_k p_k)\theta}.
$$

\begin{theorem}
For almost all $\theta$, the function $f(s, \theta)$ has analytic continuation
to the half plane $\Real s >\frac12$, where it doesn't vanish.
\end{theorem}

This theorem was proved by different authors. Here we give a proof that is based on the estimates of special
trigonometric sums.

Let $\Real s >1$. Then
$$
\ln f(s,\theta) = -\sum\limits_{p}  \ln \left(  1-\frac{e^{2\pi i\theta p}}{p^{s}}\right),
$$
$$
-\frac{f'(s,\theta)}{f(s, \theta)} = \sum_{p} \sum\limits_{m=1}^{+\infty} \frac{e^{2\pi i m\theta p} \ln p}{p^{ms}} =
$$
$$
\sum_{n=2}^{+\infty}\frac{e^{2\pi i \theta n} \Lambda(n)}{n^{s}} +
\sum_{p} \sum\limits_{m \ge 2} \frac{e^{2\pi i m\theta p} \ln p}{p^{ms}} -
\sum_{p} \sum\limits_{m\ge 2} \frac{e^{2\pi i \theta p^{m}} \ln p}{p^{ms}},
$$
where $\Lambda(n)$~--- is von Mangoldt function. Note, that the last two  series define analytic functions
for $\Real s > \frac12 $.
We show, that  for almost all $\theta$ the series
$$
F(s, \theta)= \sum_{n=2}^{+\infty}\frac{e^{2\pi i \theta n} \Lambda(n)}{n^{s}}
$$
defines analytic function on the half plane $\Real s > \frac12$. It implies that
$(\ln f(s,\theta))'$ is analytic for $\Real s >\frac12$, which yields the desired result.

To show that $F(s,\theta)$ is analytic function for $\Real s >\frac12$, it is sufficient to consider the sum
$$
S_{N} (\theta) = \sum\limits_{n=2}^{N} \Lambda(n) e^{2\pi i \theta n}
$$
and estimate it by $|S_{N} (\theta)| \ll_{\varepsilon} N^{\frac12 +\varepsilon}$ for every $\varepsilon >0$,
where the constant implied by the Vinogradov symbol $\ll_{\varepsilon}$ may depend on $\varepsilon$.
Indeed, let $s=\sigma +it$, $\sigma >\frac12 +2\varepsilon$. Then for $N>1$, Abel's partial summation formula yields
$$
\sum\limits_{n=2}^{N} \Lambda(n) e^{2\pi i \theta n} n^{-s} = -\int\limits_{1}^{N} S_{u} (\theta) d\frac{1}{u^{s}} +
S_{N} (\theta) N^{-s}
$$
$$
\ll |s| \int\limits_{1}^{N} u^{-1 -\sigma + \frac12 + \varepsilon} du +
N^{-\sigma +\frac12 + \varepsilon} < +\infty,
$$
which implies, that the series for $F(s,\theta)$ defines analytic function in the region
$\Real s >\frac12$. To complete the proof of the theorem, we prove
the following lemma.

\begin{lemma}\label{l1}
Let
$$
S_{N} (\theta) = \sum\limits_{n=2}^{N} \Lambda(n) e^{2\pi i \theta n}.
$$
Then for almost all
$\theta\in [0;1]$ and every $\varepsilon>0$ we have
$$
|S_{N} (\theta)| \ll_{\varepsilon} N^{\frac12} (\ln N)^{\frac52 + \varepsilon}.
$$
\end{lemma}

The proof of Lemma \ref{l1} is based on

\begin{lemma}\label{l2}
(Erd\"{o}s - G\'{a}l - Koksma, see \cite{Tichy}, p.154) Let $F(M,N,\theta)$ be non-negative functions
in $L^{p} [0;1]$ ($p \ge 1$) for
$M, N = 0, 1,\ldots$, $0\le \theta \le 1$ such that $F(M,0,\theta) = 0$ and
$$
F(M,N,\theta) \le F(M,N',\theta) + F(M+N', N-N',\theta)
$$
for all $M = 0,1,\ldots$ and all $0 \le N' \le N$. Suppose further
$$
\int\limits_{0}^{1} F(M,N,\theta)^p d\theta = O(g(N)),
$$
uniformly in $M = 0,1,\ldots$, where $g(N)/N$ is a non-decreasing function. Then for arbitrary
$\varepsilon >0$
$$
F(0,N,\theta)^p = O(g(N)(\ln N)^{p+1+\varepsilon})
$$
for almost all $\theta\in [0;1]$.
\end{lemma}

\begin{pfl1}
Set
$$
F(M,N, \theta)= \left| \sum\limits_{M<n\le M+N} \frac{\Lambda(n)}{\ln n} e^{2\pi i \theta n}\right| \quad
\text{if} \quad N\ge 1
$$
and
$$
F(M, 0,\theta)=0.
$$
The non-negative functions $F(M,N, \theta) \in L^2 [0;1]$ satisfy conditions of lemma \ref{l2}.
From the definition of $F(M,N, \theta)$ we find
$$
\int\limits_{0}^{1} F(M,N,\theta)^2 d\theta = \sum\limits_{M<n\le M+N} \frac{\Lambda^{2}(n)}{\ln^{2} n}\le
\sum\limits_{M<n\le M+N} 1 =N.
$$
According to lemma \ref{l2}, for every $\varepsilon>0$ and almost all $\theta\in[0;1]$ we have
$$
\left| \sum\limits_{2\le n\le N} \frac{\Lambda(n)}{\ln n} e^{2\pi i \theta n}\right| \ll_{\varepsilon}
N^{\frac12} (\ln N)^{\frac32 +\varepsilon}.
$$
By Abel's partial summation formula we find
$$
S_{N} (\theta) = \sum\limits_{2\le n\le N} \frac{\Lambda(n)}{\ln n} e^{2\pi i \theta n} \ln n =
-\int\limits_{1}^{N}  \bbC (u) d\ln u  + \bbC (N) \ln N,
$$
where
$$
\bbC (u) = \sum\limits_{1<n\le N} \frac{\Lambda(n)}{\ln n} e^{2\pi i \theta n}.
$$
Thus, for almost all $\theta\in [0;1]$ we have
$$
|S_{N} (\theta)| \ll \max\limits_{1\le u \le N} |\bbC (u)| \ln N \ll N^{\frac12} (\ln N)^{\frac52 + \varepsilon}.
$$
Obviously, once the estimate of Lemma 1 is proved for every $\varepsilon>0$ and almost all $\theta$,
it is also valid for almost all $\theta$ and every $\varepsilon>0$.
This completes the proof of lemma \ref{l1}.
\end{pfl1}


\end{document}